\documentclass[10pt,a4paper]{amsart}
\usepackage{graphicx}
\usepackage{mathrsfs}
\usepackage{color}
\usepackage[latin1]{inputenc}

\usepackage{amsmath,amssymb, amsthm, latexsym}
\usepackage{hyperref}
\input xy
\xyoption{all}

\author{S\'ebastien Alvarez}
\date{}
\title{Discretization of harmonic measures for foliated bundles}

\begin{document}

\newcounter{theorem}[section]
\newtheorem{exemple}{\bf Exemple \rm}
\newtheorem{exercice}{\bf Exercice \rm}
\newtheorem{conj}[theorem]{\bf Conjecture}
\newtheorem{defi}[theorem]{\bf Definition}
\newtheorem{lemma}[theorem]{\bf Lemma}
\newtheorem{proposition}[theorem]{\bf Proposition}
\newtheorem{coro}[theorem]{\bf Corollary}
\newtheorem{theorem}[theorem]{\bf Theorem}
\newtheorem*{mtheorem}{\bf Main Theorem}
\newtheorem*{ptheorem}{\bf Théorème Principal}
\newtheorem*{mcoro}{\bf Corollary}
\newtheorem*{pcoro}{\bf Corollaire}
\newtheorem{rem}[theorem]{\bf Remark}
\newtheorem{ques}[theorem]{\bf Question}
\newtheorem{propr}[theorem]{\bf Property}
\newtheorem{question}{\bf Question}
\def\bp{\noindent{\it Proof. }}
\def\ep{\noindent{\hfill $\fbox{\,}$}\medskip\newline}
\renewcommand{\theequation}{\arabic{section}.\arabic{equation}}
\renewcommand{\thetheorem}{\arabic{section}.\arabic{theorem}}
\newcommand{\eps}{\varepsilon}
\newcommand{\disp}[1]{\displaystyle{\mathstrut#1}}
\newcommand{\fra}[2]{\displaystyle\frac{\mathstrut#1}{\mathstrut#2}}
\newcommand{\dif}{{\rm Diff}}
\newcommand{\homeo}{{\rm Homeo}}
\newcommand{\Per}{{\rm Per}}
\newcommand{\Fix}{{\rm Fix}}
\newcommand{\A}{\mathcal A}
\newcommand{\Z}{\mathbb Z}
\newcommand{\Q}{\mathbb Q}
\newcommand{\R}{\mathbb R}
\newcommand{\C}{\mathbb C}
\newcommand{\N}{\mathbb N}
\newcommand{\T}{\mathbb T}
\newcommand{\U}{\mathbb U}
\newcommand{\D}{\mathbb D}
\newcommand{\PP}{\mathbb P}
\newcommand{\Sp}{\mathbb S}
\newcommand{\K}{\mathbb K}
\newcommand{\car}{\mathbf 1}
\newcommand{\g}{\mathfrak g}
\newcommand{\gs}{\mathfrak s}
\newcommand{\h}{\mathfrak h}
\newcommand{\rr}{\mathfrak r}
\newcommand{\s}{\sigma}
\newcommand{\fhi}{\varphi}
\newcommand{\ffhi}{\tilde{\varphi}}
\newcommand{\moins}{\setminus}
\newcommand{\ds}{\subset}
\newcommand{\W}{\mathcal W}
\newcommand{\WW}{\widetilde{W}}
\newcommand{\F}{\mathcal F}
\newcommand{\G}{\mathcal G}
\newcommand{\CC}{\mathcal C}
\newcommand{\RR}{\mathcal R}
\newcommand{\DD}{\mathcal D}
\newcommand{\M}{\mathcal M}
\newcommand{\B}{\mathcal B}
\newcommand{\HH}{\mathcal H}
\newcommand{\Hyp}{\mathbb H ^2}
\newcommand{\UU}{\mathcal U}
\newcommand{\Pp}{\mathcal P}
\newcommand{\QQ}{\mathcal Q}
\newcommand{\E}{\mathcal E}
\newcommand{\GG}{\Gamma}
\newcommand{\LL}{\mathcal L}
\newcommand{\diam}{{\rm diam}}
\newcommand{\diag}{{\rm diag}}
\newcommand{\Jac}{{\rm Jac}}
\newcommand{\Plong}{{\rm Plong}}
\newcommand{\Tr}{{\rm Tr}}
\newcommand{\Conv}{{\rm Conv}}
\newcommand{\Ext}{{\rm Ext}}
\newcommand{\Spec}{{\rm Sp}}
\newcommand{\Supp}{{\rm Supp}\,}
\newcommand{\Grass}{{\rm Grass}}
\newcommand{\Ad}{{\rm Ad}}
\newcommand{\ad}{{\rm ad}}
\newcommand{\Aut}{{\rm Aut}}
\newcommand{\Leb}{{\rm Leb}}
\newcommand{\End}{{\rm End}}
\newcommand{\cc}{{\rm cc}}
\newcommand{\grad}{{\rm grad}}
\newcommand{\proj}{{\rm proj}}
\newcommand{\mass}{{\rm mass}}
\newcommand{\dive}{{\rm div}}
\newcommand{\dist}{{\rm dist}}
\newcommand{\im}{{\rm Im}}
\newcommand{\re}{{\rm Re}}
\newcommand{\codim}{{\rm codim}}
\newcommand{\Map}{\longmapsto}
\newcommand{\vide}{\emptyset}
\newcommand{\tr}{\pitchfork}
\newcommand{\ssl}{\mathfrak{sl}}

\newenvironment{demo}{\textbf{Proof.}}{\quad \hfill $\square$}
\newenvironment{pdemo}{\textbf{Proof of the proposition.}}{\quad \hfill $\square$}
\newenvironment{IDdemo}{\textbf{Id\'ee de preuve.}}{\quad \hfill $\square$}

\def\to{\mathop{\rightarrow}}
\def\act{\mathop{\curvearrowright}}
\def\To{\mathop{\longrightarrow}}
\def\Sup{\mathop{\rm Sup}}
\def\Max{\mathop{\rm Max}}
\def\Inf{\mathop{\rm Inf}}
\def\Min{\mathop{\rm Min}}
\def\lims{\mathop{\overline{\rm lim}}}
\def\limi{\mathop{\underline{\rm lim}}}
\def\egal{\mathop{=}}
\def\dans{\mathop{\subset}}
\def\surj{\mathop{\twoheadrightarrow}}

\begin{abstract}
\textbf{Discrétisation des mesures harmoniques pour des fibrés feuilletés}

\textbf{Systèmes dynamiques}

On prouve dans cette note qu'il y a, pour certains fibrés feuilletés, une correspondance bijective entre les mesures harmoniques au sens de Garnett et les mesures sur la fibre qui sont stationnaires pour une certaine mesure de probabilité sur le groupe d'holonomie. Nous en déduisons l'unicité de la mesure harmonique pour certains feuilletages transverses à une fibration projective.
\newline

We prove in this note that there is, for some foliated bundles, a bijective correspondance between Garnett's harmonic measures and measures on the fiber that are stationary for some probability measure on the holonomy group. As a consequence, we show the uniqueness of the harmonic measure in the case of some foliations transverse to projective fiber bundles
\end{abstract}

\maketitle

\section{Version française abrégée}
Le but de cette note est de donner une application directe de la discrétisation du mouvement Brownien de Furstenberg-Lyons-Sullivan (\cite{Fu1}, \cite{LS}) à l'étude des mesures harmoniques des feuilletages introduites par Garnett \cite{Gar}. La motivation principale de ce travail est d'obtenir des théorèmes d'unicité de mesures harmoniques pour des fibrés feuilletés transversalement projectifs, généralisant ainsi certains résultats de Bonatti et G\'omez-Mont \cite{BG-M}.  Ces mesures harmoniques sont, par définition, invariantes en moyenne par diffusion Brownienne le long des feuilles. Lorsque le feuilletage est transverse à une fibration il est possible de désintégrer une telle mesure dans les fibres par rapport au volume dans la base, et de regarder l'action du groupe d'holonomie sur les mesures conditionnelles. 

\'Etant donnée une mesure de probabilité sur le groupe d'holonomie, on peut s'intéresser aux mesures sur une fibre qui sont invariantes en moyenne par l'action du groupe. Ce sont ces mesures que l'on appelle stationnaires. Nous savons, depuis les travaux de Furstenberg \cite{Fu2} qu'il existe des critères explicites garantissant l'unicité de mesures stationnaires pour des actions projectives. L'idée est alors de ramener l'étude des mesures harmoniques du feuilletage à celle des mesures stationnaires pour une certaine mesure de probabilité sur le groupe d'holonomie qui est obtenue en discrétisant le mouvement Brownien.

Le résultat principal de cette note est qu'il y a une bijection entre ces deux types de mesures. \'Enonçons précisément le théorème principal. Considérons $S$ une variété riemannienne de volume fini dont la géométrie est bornée. Soit $\Pi: N\to S$ un fibré localement trivial dont la fibre $V$ est une variété différentielle compacte. Supposons de plus que le fibré soit plat: il existe un feuilletage $\F$ de $N$ dont les feuilles sont transverses à la fibration et revêtent la base, et dont l'holonomie est donnée par une action par difféomorphismes du groupe fondamental de $S$ sur $V$.

\begin{ptheorem}
Il existe une mesure de probabilité $\mu$ sur $\GG=\pi_1(S)$ telle que pour tout fibré feuilleté $\Pi: N\to S$ dont la fibre $V$ est compacte, on ait une bijection entre les mesures harmoniques du feuilletage $\F$, et les mesures sur $V$ qui sont $\mu$-stationnaires.
\end{ptheorem}

Appliquons ce théorème au cas où la fibre est un espace projectif, et où la base est une variété compacte de courbure négative. Supposons de plus que l'action sur la fibre soit projective, contractante et fortement irréductible (les définitions précises se trouvent dans l'introduction). On obtient alors le résultat suivant qui est un corollaire du théorème principal, ainsi que d'un résultat de Guivarc'h et Raugi \cite{GR}:

\begin{pcoro}
Soit $S$ une variété riemannienne compacte à courbure sectionnelle négative et $V=\C\PP^{d-1}$, où $d\geq 2$. Considérons un fibré feuilleté $\Pi: N\to S$ de fibre $V$. Supposons que l'action de $\GG=\pi_1(S)$ sur $V$ donnée par l'holonomie du feuilletage soit contractante et fortement irréductible. Alors le feuilletage $\F$ possède une unique mesure harmonique.
\end{pcoro}

\section{Introduction}

\subsection{Foliated bundles}

Let $M$ be a $C^{3}$ complete connected and simply connected riemannian manifold with bounded geometry. We consider $\GG$ a discrete subgroup of isometries of $M$ with finite covolume, so that the quotient $S=M/\GG$ is a manifold of finite volume, and we have the identification $\GG=\pi_1(S)$. Then it is known that harmonic functions of $S$ are constant and the Brownian motion on $S$ is recurrent.

Assume that $\GG$ acts by diffeomorphisms on a compact differential manifold $V$: there is a representation $\rho:\Gamma\to\dif(V)$ (note that the differential structure on $V$ plays no role: we could have taken any compact topological space). Then $\GG$ acts diagonally on the product $M\times V$. It means that for a couple $(p,t)\in M\times V$, we define $\gamma(p,t)=(\gamma p,\rho(\gamma)t)$. It is proved in \cite{CL} that the quotient $N$ is a differential manifold, called the \textbf{suspension} of the action, and is endowed with:

\begin{itemize}
\item a structure of fiber bundle over $S$ with $V$-fibers,
\item a transversal foliation $\F$, called the \textbf{suspension foliation}, whose leaves are covering spaces of the base $S$, and whose holonomy group is given by $\widehat{\rho}(\gamma)=\rho(\gamma^{-1})$.
\end{itemize}

We can lift the metric of $S$ to the leaves of $\F$ via the fibration: this leafwise metric varies continuously with the leaf in the $C^3$-topology. Therefore, it is possible to define a Laplace operator on each leaf, that gives us a foliated Laplace operator on denoted by $\Delta_{\F}$. More precisely, if $f$ is a continuous function defined on $N$ that is $C^2$ along the leaves, $\Delta_{\F}f$ is the continuous function whose restriction to any leaf $L$ is the Laplacian of $f_{|L}$. That leads us to consider leafwise heat equation, and finally Brownian diffusion: see \cite{Gar}. The existence of probability measures on $N$ invariant by this process, called \textbf{harmonic measures} for $\F$, is guaranteed by a theorem of L.Garnett \cite{Gar} in the case where the base $S$ is compact. In the non compact case, this is a consequence of the main theorem, because, as mentioned below, stationary measures on compact spaces always exist. They are, by definition, the measures that vanish on all the Laplacians.

If $m$ is a harmonic measure for $\F$, its projection on $S$ is also harmonic (it vanishes on all the  Laplacians). But since $S$ has finite volume and bounded geometry, there exists, up to a scalar factor, only one harmonic measure: the \emph{normalized} Lebesgue measure denoted by $\Leb$. Hence, it is possible to disintegrate $m$ on the fibers with respect to Lebesgue measure in the sense of Rokhlin (see \cite{Ro}): there exists, for Lebesgue-almost every $z\in S$ a probability measure $m_z$, called the conditional measure, on the fiber $V_z\simeq V$ such that for any Borelian $B\dans N$, we have:
$$m(B)=\int_S m_z(B\cap V_z) d\Leb(z).$$

Now assume that there is a probability measure $\mu$ on $\GG$. A probability measure $\nu$ on $V$ is said to be $\mu$-\textbf{stationary} if the action of $\GG$ leaves $\nu$ invariant \emph{on average}: $\nu=\sum_{\gamma} \mu(\gamma)\gamma\ast\nu$. Existence of such measures is guaranteed by a Krylov-Bogoliubov type argument because $V$ is compact. 

While the conditional measures of harmonic probabilities are invariant on average by holonomy transportation along Brownian paths, stationary measures are invariant on average under the action of the holonomy group. The theorem that follows is our main result, and states that these objects are the same, provided a good choice of a probability measure on $\GG$.

\begin{mtheorem}
There exists a probability measure $\mu$ on $\GG$ such that for any such foliated bundle, there is a bijective correspondance between harmonic measures for the foliation $\F$ and $\mu$-stationary measures on $V$.
\end{mtheorem}

Let's consider a particular case. Assume that $S$ is a compact manifold with negative curvature, that $V=\C\PP^{d-1}$ with $d\geq 2$, and that we have a projective representation $\rho:\GG\to PSL_d(\C)$.

\begin{defi}
\begin{enumerate}
\item We say that the representation $\rho$ is \textbf{contracting} if for any probability measure on $\C\PP^{d-1}$, there exists a sequence $(\gamma_n)_{n\in\N}\in\GG^{\N}$ such that $\rho(\gamma_n)\ast m$ converges to a Dirac mass for the weak-$*$ topology.
\item We say that the representation $\rho$ is \textbf{strongly irreducible} if no finite family $\{V_1,...,V_k\}$ of proper projective subspaces of $\C\PP^{d-1}$ is invariant under every $\rho(\gamma)$, $\gamma\in\GG$.
\end{enumerate}
\end{defi}

We can state the following theorem (slightly modified to fit to our context) of Guivarc'h and Raugi \cite{GR}:

\begin{theorem}[Guivarc'h-Raugi]
\label{gura}
Assume that $\mu$ is a probability measure on $\GG$ with full support. Let's consider a contracting and irreducible representation $\rho:\GG\to PSL_d(\C)$. Assume moreover that $\rho$ satisfies the following integrability condition:
$$\sum_{\gamma\in\GG}\mu(\gamma)\log||\rho(\gamma)||<\infty.$$
Then there is a unique $\mu$-stationary probability measure $\nu$ on $\C\PP^{d-1}$.
\end{theorem}

As we will see in later, the measure $\mu$ on $\Gamma$ given by the main theorem satisfies the full support condition, aswell as the integrability condition \emph{for any} representation. Therefore, the main theorem and the theorem \ref{gura} give the following corollary:

\begin{mcoro}
Assume that $S$ is a compact riemannian manifold with negative sectional curvature. Let $\rho:\pi_1(S)\to PSL_d(\C)$ be a contracting and strongly irreducible projective representation. Let $\Pi:N\to S$ be the foliated bundle obtained by the suspension of $\rho$. Then there exists a unique harmonic measure for the suspension foliation.
\end{mcoro}

Note that, at least when $d=2$ and $S$ is a compact surface, the condition of being contracting and strongly irreducible is open and dense in the variety of representations of $\pi_1(S)$ in $PSL_2(\C)$ (see \cite{BG-MV}). Even if it seems plausible that this condition remains open and dense in more generality, we did not find any reference in the litterature.

\subsection{Discretization of the Brownian motion}

Remember that $M$ is a $C^{3}$ complete connected and simply connected riemannian manifold with bounded geometry, and  $\GG<Isom^+(M)$ is discrete of finite covolume. By discretization of Brownian motion on $M$, we mean discretization of functions that are invariant under Brownian diffusion: the harmonic functions. We fix $p_0$ once for all.

Before stating the theorem, we need two notations.
\begin{itemize}
\item The set $\HH^+(M,\Delta)$ consists of positive harmonic functions on $M$ for the Laplace-Beltrami operator $\Delta$.
\item If $(\mu_{\gamma p_0})_{\gamma\in \GG}$ is a family of probabilities on $\GG$, we say a function $h:\GG p_0\to\R^+$ is $(\mu_{\gamma p_0})_{\gamma\in\GG}$-\emph{harmonic} if for every $\gamma\in \GG$,
$$h(\gamma p_0)=\sum_{\xi\in\GG} \mu_{\gamma p_0}(\xi)h(\xi p_0).$$
\end{itemize}
We call $\HH^+(\GG p_0,(\mu_{\gamma p_0})_{\gamma\in \GG})$ the set of such functions.

\begin{theorem}[Furstenberg-Lyons-Sullivan] 
\label{discr}
There exists a family of probablity measures on $\Gamma$ denoted by $(\mu_p)_{p\in M}$ such that:
\begin{enumerate}
\item Each $\mu_p$ has full support: for any $p\in M$ and $\gamma\in\GG$, we have $\mu_p(\gamma)>0$.
\item The familly is $\GG$-equivariant, that is for every $p\in M$ and $\gamma_1,\gamma_2\in\GG$, $\mu_{\gamma_1 p}(\gamma_1 \gamma_2)=\mu_p(\gamma_2)$,
\item For any function $h\in\HH^+(\Gamma p_0,(\mu_{\gamma p_0})_{\gamma\in \GG})$, the formula 
$$\Phi h(p)=\sum_{\gamma\in\GG} \mu_p(\gamma) h(\gamma p_0)$$
 defines a smooth harmonic function of $p\in M$.
\item
The application $\Phi:\HH^+(\GG p_0,(\mu_{\gamma p_0})_{\gamma\in \GG})\to\HH^+(M,\Delta)$ is a bijection, the inverse application just being the restriction.
\end{enumerate}
\end{theorem}

In their article \cite{BL}, Ballmann and Ledrappier proved that for spaces with sectional curvature pinched into two constants, the family $(\mu_p)_{p\in M}$ can be chosen with the following  property of integrability:

\begin{theorem}[Ballmann-Ledrappier]
\label{tek}
If the sectional curvature of $M$ is pinched into two negative constants, it is possible to choose the family $(\mu_p)_{p\in M}$ of the theorem \ref{discr} such that each measure $\mu_p$ on $\GG$ has finite first moment, that is:
$$\sum_{\gamma\in\GG} \mu_p(\gamma)\dist(p_0,\gamma p_0)<\infty.$$
\end{theorem}

\section{Proof of the theorem}

In what follows, we are under the hypothesis of the main theorem: $M$ is a $C^3$ complete connected and simply connected riemannian manifold with bounded geometry, $\GG$ is a lattice of direct isometries of $M$ so that $S=M/\GG$ is a manifold with finite volume. $V$ is a compact differential manifold and $\Pi:N\to S$ is a foliated bundle obtained by suspension of a representation $\rho:\GG\to\dif(V)$. For sake of simplicity, if $\nu$ is a measure on $V$, we will use the notation $\gamma\ast\nu=\rho(\gamma)\ast\nu$. We have fixed once for all a point $p_0\in M$, and its projection $z_0\in S$. Finally, $(\mu_p)_{p\in M}$ is the familly given by the theorems of the previous section.

\subsection{Conditional measures of harmonic measures are stationary}

The condition for having a foliation transversal to the bundle $\Pi:N\to S$ is equivalent to that of flatness. This means that there exists a locally finite family of small balls of $S$, say $(U_i)_{i\in\N}$, that trivialize the bundle: $\Pi^{-1}(U_i)\simeq U_i\times V$, such that the intersection of any two of the $U_i$'s is connected, and the transition functions are of the form $(z,t)\in (U_i\cap U_j)\times V\mapsto(z,\tau_{ij}(t))\in (U_i\cap U_j)\times V$, $\tau_{ij}$ being a diffeomorphism of $V$ independent of $z$.

These $\tau_{ij}$ are called the \textbf{holonomy transformations}. We have the following relations: $\tau_{ji}=\tau_{ij}^{-1}$ and if $U_i\cap U_j\cap U_k\neq\vide$, then $\tau_{jk}\tau_{ij}=\tau_{ik}$. If we have a continuous path $c$ in the base, and a chain  $(U_{i_0},...,U_{i_n})$ of discs that covers $c$, the holonomy transportation along $c$ is $\tau_c=\tau_{i_{n-1}i_n}...\tau_{i_0i_1}$. This application does not depend on the covering, nor on the choice of $c$, but only on the homotopy class of $c$. Note that if $\gamma$ is an element of $\GG$ identified with a homotopy class of a loop based in $z_0$, then $\tau_{\gamma}=\rho(\gamma)^{-1}$. In the sequel, we will use the abusive notation $\Pi^{-1}(U_i)=U_i\times V$. Note that it is possible to assume that the $U_i$'s are sufficiently small so that they also trivialize the universal cover $M\to S$. 

Let $m$ be a harmonic measure for the foliation and $(m_z)_{z\in S}$ its disintegration with respect to Lebesgue: each $m_z$ is a probability measure on the fiber $V_z$. Remember that a theorem of Rokhlin \cite{Ro} gives uniqueness of this familly up to a zero measure of $z\in S$. It is possible to show (cf. \cite{Gar}) that in a chart $U_i\times V$, we have the following decomposition of the measure:
\begin{equation}
\label{boite}
m_{|U_i\times V}=h_i(z,t) \Leb(z)\, \nu_i(t),
\end{equation}
where $\nu_i$ is a measure on $V$, and $h_i$ is a measurable function, defined on $U_i\times A_i$ where $A_i\dans V$ has full $\nu_i$-measure, which is of class $C^2$ and harmonic on almost-each plaque. Hence, it allows us to indentify the conditional measures in the following way: for $z\in U_i$,
$$m_z=h_i(z,t)\,\nu_i(t).$$

Note that if $c$ is a continuous path on the basis starting at $z$, then the holonomy acts on the measures $m_z$ by $\tau_c\ast m_z=h(z,\tau_c^{-1}(t))\,\tau_c\ast\nu_i(t)$. The goal of this section is to prove the following:

\begin{proposition}
\label{stat}
Let $z_0$ be the projection of $p_0$ on $S$. Then the measure $m_{z_0}$ is $\mu_{p_0}$-stationary:
$$m_{z_0}=\sum_{\gamma\in\GG} \mu_{p_0}(\gamma)\,\,\gamma\ast m_{z_0}.$$
\end{proposition}

In order to prove this proposition, we have to state explicitely the cocycle relations introduced by the holonomy. If we evaluate $m$ over $U_i\cap U_j$, we obtain that for $t\in A_i\cap \tau_{ij}^{-1}(A_j)$:

\begin{equation}
\label{coc1}
\frac{h_i(z,t)}{h_j(z,\tau_{ij}(t))}=\frac{d[(\tau_{ij}^{-1})\ast\nu_j]}{d\nu_i}(t).
\end{equation}
This shows that even if the $\nu_i$ are not invariant under holonomy, the \emph{class} of the measure is. Therefore, it is possible to assume that $\tau_{ij}(A_i)=A_j$. This is the same as taking $A\dans V$ full for each measure $\nu_i$ and invariant under holonomy maps, in such a way that each $h_i$ is defined on $U_i\times A$.

We have chosen $U_0$ small enough so that all the plaques $U_0\times\{t\}$ have a section in the universal cover $\widetilde{U_0}\dans M$ that contains $p_0$: it is possible to lift the function $h_0(.,t)$ in order to obtain a harmonic function $H_t:\widetilde{U_0}\to\R$. The next lemma shows that this function can be extended harmonically to the whole $M$.
\begin{lemma}
For any $t\in A$, the function $H_t$ can be extended harmonically to $M$. Moreover, if $p\in M$ and $c$ is the projection on $S$ of the segment $[p_0,p]$, and if $z\in U_i$ is the ending point of $c$, then we have the relation:
\begin{equation}
\label{cocycle}
H_t(p)=\frac{d[(\tau_c^{-1})\ast\nu_i]}{d\nu_0}(t) h_i(z,\tau_c(t)).
\end{equation}
\end{lemma}

\begin{demo}
We have to show by induction on the number of different $U_i$'s met by the projection $c$ of $[p_0,p]$. If this number is zero, there is nothing to prove since the function $H_t$ is harmonic on $\widetilde{U_0}$. Now for the heredity, it is enough to remark that, by the relation \ref{coc1} if $U_i\cap U_j\neq\vide$, then the function $h_{ij}:z\in U_j\times\{\tau_{ij}(t)\}\mapsto{h_j(z,\tau_{ij}(t))}\frac{d[(\tau_{ij}^{-1})\ast\nu_j]}{d\nu_i}(t)$ is a harmonic continuation of $h_i(.,t):U_i\times\{t\}\to\R$.
\end{demo}

We will need the next lemma which is a consequence of the previous one and from which we will obtain the proof of the proposition, aswell as the fact that the map that associates to a harmonic measure $m$ its conditional measure on the fiber $V_{z_0}$ is injective.

\begin{lemma}
\label{mainlemma}
Take $p\in M$, and consider the projection $c$ on $S$ of the segment $[p_0,p]$. Then for all $t\in A$:
\begin{equation}
\label{discret}
H_t(p)=\sum_{\gamma\in\GG}\mu_p(\gamma)h_0(z_0,\rho(\gamma)^{-1} t)\frac{d\gamma\ast\nu_0}{d\nu_0}(t)
\end{equation}
\end{lemma}

\begin{demo}
Since for any $t\in A$, the function $H_t$ is harmonic in $M$, it is possible to use the theorem \ref{discr}: for every $t\in A$ and $p\in M$, we have $$H_t(p)=\sum_{\gamma\in\GG}\mu_p(\gamma)H_t(\gamma p_0).$$
The homotopy class of the projection of $[p_0,\gamma p_0]$ is of course given by $\gamma$, and the associated holonomy transformation is given by $\tau_{\gamma}=\rho(\gamma)^{-1}$. Hence, by the relation \ref{cocycle}, $H_t(\gamma p_0)=h_0(z_0,\rho(\gamma)^{-1}t)\frac{d\gamma\ast\nu_0}{d\nu_0}(t) $. Finally, the relation \ref{discret} is true.
\end{demo}

Now, we can end the proof of the proposition \ref{stat}:

\begin{pdemo}
Since $A$ is full for $\nu_0$, we can multiply the relation \ref{discret} by the measure $\nu_0$. When $p=p_0$, the left member becomes $H_t(p_0)\nu_{0}$, that is by definition, $h_0(z_0,t)\nu_0$. But remember that the latter is equal to $m_{z_0}$.

The right member is a combination of the measures $h_0(z_0,\rho(\gamma)^{-1}t)\,\frac{d\gamma\ast\nu_0}{d\nu_0}(t)\nu_0 =h_0(z_0,\rho(\gamma)^{-1}t)\,\gamma\ast\nu_0$ with weight $\mu_{p_0}(\gamma)$. But remember also that the latter measure is $\gamma\ast m_{z_0}$, we obtain:
$$m_{z_0}=\sum_{\gamma\in\GG} \mu_{p_0}(\gamma)\,\,\gamma\ast m_{z_0},$$
concluding the proof of the proposition.
\end{pdemo}

Similarly, if we apply the lemma \ref{mainlemma} for any $p$:

\begin{lemma}
\label{conditional}
Take $p\in M$, and call $c$ the projection on $S$ of the segment $[p_0,p]$, and $z$ the projection of $p$ on $S$. Then the following equality is true:
$$m_z=\tau_c\ast\left(\sum_{\gamma\in\GG}\mu_p(\gamma) \gamma\ast m_{z_0}\right)$$
\end{lemma}

\begin{demo}
Let's $p$,$c$ and $z$ be as in the lemma. We assume that $z\in U_i$. For all $t\in A$, the relation \ref{cocycle} gives $H_t(p)=\frac{d(\tau_c^{-1}\ast\nu_i)}{d\nu_0}(t) h_i(z,\tau_c(t))$. Hence, if we multiply the relation \ref{discret} by $\nu_0$, the left member becomes $h_i(z,\tau_c(t)) \frac{d(\tau_c^{-1}\ast\nu_i)}{d\nu_0}(t)\,\nu_0= h_i(z,\tau_c(t))\,\tau_c^{-1}\ast\nu_i$. But remember that the latter is equal to $\tau_c^{-1}\ast m_z$.

On the other hand, the second member is transformed exactly in the same way as in the proof of the proposition, the only difference is that we obtain the combination of the $\gamma\ast m_{z_0}$ whose weights are given by the $\mu_{p}(\gamma)$ and no longer by the $\mu_{p_0}(\gamma)$. The relation we obtain is:
$$\tau_c^{-1}\ast m_z=\sum_{\gamma\in\GG}\mu_p(\gamma) \gamma\ast m_{z_0}.$$
We conclude the proof of the proposition by pushing by $\tau_c$.
\end{demo}

Hence, all conditional measures of $m$ can be constructed from the one on the fiber of $z_0$. We therefore have the folowing injectivity result:
\begin{coro}
\label{injectivity}
Let's assume the hypothesis of the proposition. Let $m$ and $m'$ be two harmonic measures such that the conditional measures $m_{z_0}$ and $m_{z_0}'$ are equal. Then the two measures $m$ and $m'$ are equal.
\end{coro}

\subsection{Construction of a harmonic measure from a stationary measure}

Let $\nu$ be a $\mu_{p_0}$-stationary measure on $V$. Remember that this implies in particular that all the $\gamma\ast\nu$ are in the same measure class. In order to construct a harmonic measure for $\F$, we will construct a measure on each $\{p\}\times V$, $p\in M$, thanks to the family of theorem \ref{discr}, integrate it with respect to the Lebesgue measure, and pass to the quotient. The lemma \ref{conditional} gives us candidates for conditional measures. Indeed, if $p\in M$, we state:

\begin{gather}
\label{def}
\tag{$\star$}
\tilde{m}_p=\sum_{\gamma\in\GG}\mu_p(\gamma)\gamma\ast \nu.
\end{gather}
Note that in particular, by stationarity of the measure $\nu$, we have $\tilde{m}_{p_0}=\nu$.

\begin{proposition}
Let $\tilde{m}$ be a measure on $M\times V$ obtained as the integration of the $\tilde{m}_p$ with respect to Lebesgue. Then $\tilde{m}$ is harmonic and passes to the quotient by the diagonal action, giving a harmonic measure $m$  for the suspension $\F$ on $N$.
\end{proposition}

\begin{demo}
Let's prove first the second part of the proposition: we have to show that $\tilde{m}$ is invariant under the action of $\GG$. First, note that since $\GG$ acts on $M$ by isometries, its action on $M$ leaves invariant the Lebesgue measure. Hence, we are left to show that the conditional measures are preserved by the action of $\GG$, we want to show that for any $\xi\in\Gamma$, 
\begin{equation}
\label{invariance}
\xi\ast\tilde{m}_{p}=\tilde{m}_{\xi p}.
\end{equation}
And this is a simple consequence of the equivariance (see item 2 of the theorem \ref{discr}) of the family of measures $(\mu_p)_{p\in M}$. For any $p\in M$ and $\xi\in\GG$,

$$\xi\ast\tilde{m}_p=\sum_{\gamma\in\GG} \mu_{p}(\gamma)\,\,(\xi\gamma)\ast\nu
                  =\sum_{\gamma\in\GG} \mu_{\xi p}(\xi\gamma)\,\,(\xi\gamma)\ast\nu
                  =\sum_{\eta\in\GG} \mu_{\xi p}(\eta)\,\,\eta\ast\nu
                  =\tilde{m}_{\xi p}.$$

We then have to prove the first part. It is sufficient to show that in restriction to any $M\times\{t\}$ the measure $\tilde{m}$ has a harmonic density with respect to the Lebesgue measure. Hence, let $A\dans V$ be a full $\nu$-measure Borel set such that for any $t\in A$ and  $\gamma\in\GG$, the derivative $\frac{d\gamma\ast\nu}{d\nu}(t)$ exists. In the next lemma we introduce a family of harmonic functions which, as we will show later, are the densities with respect to Lebesgue, of the conditional measures of $\tilde{m}$ on the $M\times\{t\}$ with $t\in A$.

\begin{lemma}
For any $t\in A$, the following function of $p\in M$ is harmonic:
$$H_t(p)=\sum_{\gamma\in\GG} \mu_p(\gamma) \frac{d\gamma\ast\nu}{d\nu}(t).$$
\end{lemma}

\begin{demo}
Remark first that by definition, for any $t\in A$ and $p\in M$, we have: $H_t(p)=\frac{d\tilde{m}_p}{d\nu}(t).$ Since by \ref{invariance} we have that for any $\xi\in\GG$, $\tilde{m}_{\xi p_0}=\xi\ast\tilde{m}_{p_0}$ and by stationarity, $\tilde{m}_{p_{0}}=\nu$, then $\tilde{m}_{\xi p_0}=\xi\ast\nu$ holds for all $\xi$. 
Hence, for any $t\in A$ and $\xi\in\GG$, $H_t(\xi p_0)=\frac{d\xi\ast\nu}{d\nu}(t)$, and by definition, for any $t\in A$ and $\gamma\in\GG$,
$$H_t(\gamma p_0)=\sum_{\xi\in\GG}\mu_{\gamma p_0}(\xi)\frac{d\tilde{m}_{\xi p_0}}{d\nu}(t)
                 =\sum_{\xi\in\GG}\mu_{\gamma p_0}(\xi)H_t(\xi p_0).$$
And this means that for any $t\in A$, the function $\gamma\in\GG\mapsto H_t(\gamma p_0)$ is harmonic for the familly $(\mu_{\gamma p_0})_{\gamma\in\GG}$ (see paragraph 2.2). By the theorem \ref{discr} of Furstenberg-Lyons-Sullivan, the function $p\mapsto H_t(p)$ is harmonic.
\end{demo}
\newline

Hence, the conditional measures of $\tilde{m}'=H_t(p) \Leb(p)\nu(t)$ on the $M\times\{t\}$ have harmonic densities with respect to Lebesgue. We can disintegrate this measure on the $\{p\}\times V$, and the conditional measures are the $\tilde{m}'_p=H_t(p)\nu(t)$, and, by construction, this measure is equal to $\tilde{m}_p$. This means that the two measures $\tilde{m}$ and $\tilde{m}'$ are equal. Therefore $\tilde{m}$ is harmonic, and the proof of the proposition is now over.
\end{demo}

Finally, if $z_0\in S$ is the projection of $p_0$, the conditional measure at $z_0$ of the harmonic measure is given by $\nu$ (it is a probability measure because $H_t(p_0)=1$ for any $t$). This fact, together with the corollary \ref{injectivity}, allows us to end the proof of the main theorem by stating the proposition:
\begin{proposition}
The application that associates to each harmonic measure for $\F$ the conditional measure at $z_0$ is a bijection from the set of harmonic measures for $\F$ to the one of measures on $V$ that are $\mu_{p_0}$-stationary.
\end{proposition}

\section{Unique ergodicity of harmonic measures}

Assume that $S$ is a $n$-dimensional compact manifold with negative sectional curvature. We consider the associated family $(\mu_p)_{p\in M}$ on $X$, and a projective representation $\rho:\GG=\pi_1(S)\to PSL_d(\C)$.

Since $S$ is compact, the two following invariant distances on $\GG$ are equivalent: $d_1(\gamma_1,\gamma_2)=\dist(p_0,\gamma_1^{-1}\gamma_2 p_0)$ and $d_2(\gamma_1,\gamma_2)$, the word distance with respect to a symetric system of generators. But for \emph{any} representation $\rho:\GG\to PSL_d(\C)$, if $||.||$ is an operator norm and $C$ is greater than every $\log||\gamma_i||$ for a symetric system of generators $(\gamma_i)$ then for any $\gamma$, $\log(||\rho(\gamma)||)\leq Cd_2(Id,\gamma)$. Hence $\log||\rho(\gamma)||=O(\dist(p_0,\gamma p_0))$, and since by theorem \ref{tek}, $\mu$ has finite moment, we have the following lemma:

\begin{lemma}
The application $\gamma \mapsto\log||\rho(\gamma)||$ is $\mu_{p_0}$-integrable.
\end{lemma}

Now, assume that the action of $\GG$ is strongly irreducible and contracting, and remember that by theorem \ref{discr}, the measure $\mu_{p_0}$ has full support. By the theorem \ref{gura}, there exists a unique $\mu_{p_0}$-stationary measure $\nu$. As a consequence of the above, we are able to prove the following theorem:

\begin{theorem}
Any foliation obtained by suspending a contracting strongly irreducible representation $\rho:\pi_1(S)\to PSL_d(\C)$, $S$ being a negatively curved compact manifold, possesses a unique harmonic measure.
\end{theorem}

Remark that we only assumed the compacity of the basis in order to have the hypothesis of integrability required by the theorem \ref{gura}. It can be dropped if we know that the representation is integrable. We claim that this is the case for example if $S$ is a negatively curved (the curvature does not need to be constant) surface with finite area and if we suppose that over each geodesic around a cusp, the holonomy map has unitary eigenvalues (we say that it is parabolic). The proof will appear in a forthcoming paper.\\

\textbf{Remerciements.}
Je tiens d'abord à remercier mon directeur de thèse Christian Bonatti pour sa relecture minutieuse des premières versions de cette note et pour ses remarques précieuses. C'est Bertrand Deroin qui m'a suggéré de revisiter les travaux de Bonatti et G\'omez-Mont sous l'angle des marches aléatoires, et François Ledrappier qui m'a expliqué la discrétisation du mouvement Brownien. Je les en remercie chaleureusement. Merci enfin au rapporteur anonyme pour sa lecture attentive.

\vspace{10pt}

\noindent \textbf{Sébastien Alvarez (sebastien.alvarez@u-bourgogne.fr)}\\
\noindent  Institut de Math\'ematiques de Bourgogne, CNRS-UMR 5584\\
\noindent  Université de Bourgogne, 21078 Dijon Cedex, France\\


\begin{thebibliography}{widest-label}

\bibitem[BL]{BL} W.Ballmann, F.Ledrappier, Discretization of positive harmonic functions on Riemannian manifolds and Martin boundary, \emph{Actes de la Table Ronde de la Géométrie Différentielle (Luminy 1992), Sémin. Congr.}, \textbf{1}, Soc. Math. Fr., Paris, (1996), 77-92.

\bibitem[BG-M]{BG-M} C.Bonatti, X.G\'omez-Mont, Sur le comportement statistique des feuilles de certains feuilletages holomorphes, \emph{Monogr. Enseign. Math.}, \textbf{38}, (2001), 15-41.

\bibitem[BG-MV]{BG-MV} C.Bonatti, X.G\'omez-Mont, M.Viana, G\'en\'ericit\'e d'exposants de Lyapunov non-nuls pour des produits d\'eterministes de matrices, \emph{Ann.I. H. Poincar\'e. AN}, \textbf{20}, (2003), 579-624.

\bibitem[CL]{CL} C.Camacho, A.Lins Neto: \emph{Geometric Theory of Foliations}, Birkhäuser, Boston Inc., 1985.

\bibitem[Fu1]{Fu1} H.Furstenberg, Random walks and discrete subgroups of Lie groups, \emph{Adv. Probab. Related Topics}, \textbf{1}, Dekker, New York (1971), 1-63.

\bibitem[Fu2]{Fu2} H.Furstenberg, Noncommuting random products, \emph{Trans. Amer. Math. Soc}, \textbf{108}, (1963), 377-428.

\bibitem[GR]{GR} Y.Guivarc'h, A.Raugi, Frontière de Furstenberg, propriété de contraction et théorèmes de convergence, \emph{Z. Wahrsch. Verw. Gebiete}, \textbf{69}, (1985), 187-242.

\bibitem[Gar]{Gar} L.Garnett, Foliations, The ergodic theorem and Brownian motion, \emph{J. Funct. Anal.}, \textbf{51}, (1983), 285-311.

\bibitem[LS]{LS} T.Lyons, D.Sullivan, Function theory, random paths and covering spaces, \emph{J. Diff. Geom.}, \textbf{19}, (1984), 299-323.

\bibitem[Ro]{Ro} V.A. Rokhlin, On the fundamental ideas of measure theory, \emph{Trans. Amer. Math. Soc.}, \textbf{10}, (1962), 1-52.


\end{thebibliography}
\end{document}